\documentclass[11pt]{article}
\usepackage{amsthm,amsmath,amssymb,anysize}
\def\qed{\hbox to 0pt{}\hfill$\rlap{$\sqcap$}\sqcup$}
\newtheorem{lemma}{Lemma}[section]
\newtheorem{theorem}[lemma]{Theorem}

\newtheorem{proposition}[lemma]{Proposition}

\newtheorem{corollary}[lemma]{Corollary}
\marginsize{3cm}{3cm}{2.3cm}{2.2cm}
\numberwithin{equation}{section}

\begin{document}
\title{\textbf{Derivations  from  the even parts into the odd
parts for
  Lie superalgebras $W$ and $S$}}
\author{Wende Liu and Baoling Guan \\
\\
 {Department of Mathematics,} \\
 {Harbin Normal University, Harbin 150080,   China}\\
\\
{Email addresses: wendeliu@sohu.com, guanbl-2004@163.com}}
\date{ }
\maketitle
\begin{quotation}
\small\noindent \textbf{Abstract}:   Let  $\mathcal{W}$ and
$\mathcal{S}$ denote the even parts of the general Witt
superalgebra $W$ and the special superalgebra  $S$ over a field of
characteristic $ p>3,$ respectively. In this note, using the
method of reduction on $\mathbb{Z}$-gradations, we   determine the
derivation space $\mathrm{Der}(\mathcal{W}, W_{\overline{1}})$
from $\mathcal{W}$ into  $W_{\overline{1}} $ and  the derivation
space $\mathrm{Der}(\mathcal{S}, W_{\overline{1}})$ from
$\mathcal{S}$ into  $W_{\overline{1}}. $ In particular, the
derivation space $\mathrm{Der}(\mathcal{S}, S_{\overline{1}})$ is
determined.
\\

  \\

\end{quotation}\vspace{0.5cm}

\section{Introduction}

  The underlying field $\mathbb{F}$ is assumed of characteristic $p>3 $
  throughout. We shall
 study the derivations from the even parts  of the generalized Witt superalgebra
$W$ and the special superalgebra $S$  into the odd part of $W,$
where $W_{\overline{1}}$ is viewed as modules for $ {W}_{\overline{0}}$
and $ {S}_{\overline{0}}$ by means of the
adjoint representation. The motivation came from the following
observation. Let $L=L_{\overline{0}}\oplus L_{\overline{1}} $ be
 a Lie superalgebra. Then $L_{\overline{0}}$ is a Lie algebra and
 $L_{\overline{1}}$ is an $L_{\overline{0}}$-module.
 Two questions arise naturally:  Does the derivation algebra of the even part of $L$
  coincide with the even part of the superderivation algebra of
 $L?$ Does the derivation space from $L_{\overline{0}}$ into $
 L_{\overline{1}}$ coincide with the odd part of the
 superderivation algebra of $L?$ For the generalized Witt superalgebra
 and the special superalgebra   the first question was  answered affirmatively in \cite{g6}.
 In this note,  the second question will also be answered for these two Lie superalgebras of Cartan type.
 Speaking accurately, we shall determine
  the derivation spaces from the even parts  of the generalized
  Witt superalgebra
$W$ and the special superalgebra $S$
  into the odd part of $W.$   As a direct
consequence, the derivation space from the even part into the odd part for the special superalgebra is determined.

The author would like to thank
the anonymous referee for the paper \cite{LZ3} for  posing such an interesting  question.

In this note we adopt the notation and concepts used in \cite{g6},
but here, for convenience and completeness, we repeat certain
necessary symbols and notions.

 Let $\mathbb{Z}_{2}=\{\overline{0}, \overline{1}\}$ be
the field of two elements. For a   vector superspace
$V=V_{\overline{0}}\oplus V_{\overline{1}}, $ we denote by
$\mathrm{p}(a)=\theta$ the \textit{parity of a homogeneous
element} $a\in V_{\theta}, \theta\in \mathbb{Z}_{2}.$
 We assume throughout  that the notation   $\mathrm{p}(x)$ implies
that
  $x$ is a $\mathbb{Z}_2$-homogeneous element.

Let $\frak{g}$ be a  Lie algebra and $V$   a $\frak{g}$-module. A
linear mapping $D:\frak{g}\rightarrow V$ is called a
\textit{derivation } from $\frak{g}$  into $V$ if $D(xy)=x\cdot
D(y)-y\cdot D(x)$ for all $x,y\in \frak{g}.$ A derivation
$D:\frak{g}\rightarrow V$ is called \textit{inner} if there is
$v\in V $ such that $D(x)=x\cdot v$ for all $x\in \frak{g}.$
Following \cite[p. 13]{SF}, denote by $\mathrm{Der}(\frak{g},V)$
the \textit{derivation space} from $\frak{g}$  into $V.$  Then
$\mathrm{Der}(\frak{g},V)$ is a $\frak{g}$-submodule of
$\mathrm{Hom}_{\mathbb{F}}(\frak{g},V).$ Assume in addition that
$\frak{g}$ and $V$ are finite-dimensional and that
$\frak{g}=\oplus_{r\in \mathbb{Z}} \frak{g}_{r}$ is
$\mathbb{Z}$-graded and $V=\oplus_{r\in\mathbb{ Z}} V_{r}$ is a
$\mathbb{Z}$-graded $\frak{g}$-module. Then
$\mathrm{Der}(\frak{g},V)=\oplus_{r\in \mathbb{Z}}
\mathrm{Der}_{r}(\frak{g},V)$ is a $\mathbb{Z}$-graded
$\frak{g}$-module  by setting
$$\mathrm{Der}_{r}(\frak{g},V):=\{D\in
\mathrm{Der}(\frak{g},V)\mid D(\frak{g}_{i})\subset V _{r+i} \
\mbox{for all} \ i \in\mathbb{Z}\}.$$ In the case of $V=\frak{g},$
the \textit{derivation algebra} $\mathrm{Der}(\frak{g})$ coincides
with $\mathrm{Der}(\frak{g},\frak{g}) $ and
$\mathrm{Der}(\frak{g}) =\oplus _{r\in \mathbb{Z}}
 \mathrm {Der}_{r}(\frak{g})$ is a $\mathbb{Z}$-graded Lie algebra.
If $\frak g=\oplus_{-r\leq i\leq s}\frak g_i$ is a $\mathbb
    Z$-graded Lie  algebra, then $\oplus_{-r\leq i\leq 0}\frak
    g_i$ is called the \textit{top of} $\frak g$ (with respect to the
    gradation).

 In the below we review the notions of  modular Lie
superalgebras $W$ and $S$ of Cartan-type and their gradation
structures.
 In addition to the standard notation
$\mathbb {Z},$ we write  $\mathbb {N}$ for the set of positive
integers, and ${\mathbb {N}}_0$ for the set of nonnegative
integers. Henceforth, we will let  $m,n$  denote fixed integers in
$\mathbb {N}\setminus \{1,2\} $ without notice. For $\alpha = (
\alpha _1,\ldots,\alpha _m ) \in \mathbb {N}_0^m,$ we put
$|\alpha| =\sum_{i=1}^m\alpha _i.$ Let $\mathcal{O}(m)$ denote the
\textit{divided power algebra over} $\mathbb{F}$ with an ${\mathbb
F}$-basis $\{ x^{( \alpha ) }\mid \alpha \in \mathbb{N}_0^m \}. $
For $\varepsilon _i=( \delta_{i1}, \ldots,\delta _{im}),$ we
abbreviate $x^{( \varepsilon _i)}$ to $x_i,$ $i=1,\ldots,m.$ Let
$\Lambda (n)$ be the \textit{exterior superalgebra over}
$\mathbb{F}$ in $n$ variables $x_{m+1},\ldots ,x_{m+n}.$ Denote
the tensor product by $\mathcal{O}(m,n)
=\mathcal{O}(m)\otimes_{\mathbb{F}} \Lambda(n).$ Obviously,
$\mathcal{O}(m,n)$ is an associative superalgebra with a
$\mathbb{Z}_2$-gradation induced by the trivial
$\mathbb{Z}_2$-gradation of $\mathcal{O}(m)$ and the natural
$\mathbb{Z}_2$-gradation of $\Lambda (n).$ Moreover,
$\mathcal{O}(m,n)$ is super-commutative. For $g\in \mathcal{O}(m),f\in
\Lambda(n),$ we write $gf $  for $ g\otimes f.$   The following
formulas hold in $\mathcal{O}(m,n):$
$$
x^{(\alpha) }x^{(\beta) }=\binom{\alpha +\beta }{ \alpha}
 x^{( \alpha +\beta)}\quad\mbox{for}\  \alpha,\beta \in
{\mathbb N}_0^m;
$$
$$
x_kx_l=-x_lx_k\quad\mbox{for}\ k,l=m+1,\ldots,m+n;
$$
$$
x^{( \alpha ) }x_k=x_kx^{( \alpha ) }\quad\mbox{for}\ \alpha \in
\mathbb{N}_0^m, k=m+1,\ldots,m+n,
$$
where $\binom{ \alpha +\beta} {\alpha}:=\prod_{i=1}^m\binom{
\alpha _i+\beta _i}{ \alpha _i}.$ Put $Y_0:=\{ 1,2,\ldots,m \},$
$Y_1:=\left\{m+1,\ldots,m+n\right\} $ and $Y:=Y_0\cup Y_1.$ For
convenience, we adopt the notation $r':=r+m$ for $r\in Y_{1}.$
Thus, $Y_1:=\left\{ 1', 2', \ldots,n'\right\}.$ Set
$$\mathbb{B}_k:=\left\{ \langle i_1,i_2,\ldots
,i_k\rangle |m+1\leq i_1<i_2<\cdots <i_k\leq m+n\right\} $$ and $
\mathbb{B}:=\mathbb{B}(n)=\bigcup\limits_{k=0}^n\mathbb{B}_k,$
where $\mathbb{B}_0:=\emptyset.$ For $u=\langle i_1,i_2,\ldots
,i_k\rangle \in \mathbb{B}_k,$ set $ |u| :=k,| \emptyset | :=0$,
$x^\emptyset :=1,$ and $x^u:=x_{i_1}x_{i_2}\ldots x_{i_k};$
 we use  also $u$ to stand for the set $\{i_1,i_2,\ldots
,i_k\}.$ Clearly, $\left\{ x^{\left( \alpha \right) }x^u\mid\alpha
\in \mathbb{N}_0^m,u\in \mathbb{B}  \right\} $ constitutes an
$\mathbb{F}$-basis of $\mathcal{O} \left( m,n\right). $ Let
$D_1,D_2,\ldots,D_{m+n}$ be the linear transformations of
$\mathcal{O} \left( m,n\right) $ such that
\[
D_r ( x^{\left( \alpha \right) }x^u ) =\left\{
\begin{array}{l}
x^{\left( \alpha -\varepsilon _r\right) }x^u,\quad \quad \quad
\quad r\in Y_0
\\
x^{\left( \alpha \right) }\cdot \partial x^u/\partial x_r,\quad
\quad r\in Y_{1.}
\end{array}
\right.
\]
Then $D_1,D_2,\ldots,D_{m+n}$ are superderivations of the
superalgebra $\mathcal{O} \left( m,n\right).$ Let
\[
W\left(m,n\right) =\Big\{ \sum\limits_{r\in Y}f_r D_r \mid f_r\in
\mathcal{O} \left( m,n\right),r\in Y\Big\}.
\]
Then $W\left( m,n\right) $ is a Lie superalgebra, which is
contained in ${\rm Der} ( \mathcal{O} \left(m,n\right)). $ Obviously,
${\rm p}( D_i) =\tau (i),$ where
\[
\tau \left( i\right) :=\left\{
\begin{array}{l}
\overline{0},\quad \quad i\in Y_0 \\
\overline{1},\quad \quad i\in Y_{1.}
\end{array}
\right.
\]
One may verify that
\[[fD, gE]=fD(g)E-(-1)^{{\rm
p}(fD){\rm p}(gE)}gE(f)D+(-1)^{{\rm p}(D){\rm p}(g)}fg[D, E]\] for
$ f,g\in \mathcal{O} (m,n),$ $ D,E\in {\rm Der}\ \mathcal{O} (m,n ).$  Let
\[
\underline{t}:=\left( t_{1},t_{2},\ldots,t_m\right) \in
\mathbb{N}^m,\quad \pi:=\left( \pi _1,\pi _2,\ldots,\pi _m\right)
\]
where $\pi _i:=p^{t_i}-1,i\in Y_0.$ Let
$\mathbb{A}:=\mathbb{A}\left( m;\underline{t}\right) =\left\{
\alpha \in \mathbb{N}_0^m\mid\alpha _i\leq \pi _i,i\in Y_0
\right\}.$ Then
\[
\mathcal{O} \left( m,n;\underline{t}\right) :={\rm
span}_{\mathbb{F}}\left\{ x^{\left( \alpha \right) }x^u \mid
\alpha \in \mathbb{A},u\in \mathbb{B} \right\}
\]
is a finite-dimensional subalgebra of $\mathcal{O} \left(m,n\right) $
with   a natural ${\mathbb Z}$-gradation $\mathcal{O}
\left(m,n;\underline{t}\right)=\bigoplus _{r=0}^{\xi}
\mathcal{O}(m,n;\underline{t})_{r}$ by putting
        $$\mathcal{O}(m,n;\underline{t})_{r}:=
        {\rm span} _{\mathbb F}\{ x^{(\alpha)} x^{u}\mid |\alpha|+|u|=r \},
        \quad \xi:=|\pi|+n.$$
Set
\[
W\left( m,n;\underline{t}\right):=\Big\{ \sum\limits_{r\in Y}
f_rD_r \mid f_r\in \mathcal{O} \left( m,n;\underline{t}\right),r\in
Y\Big\}.
\]
Then $W\left( m,n;\underline{t}\right) $ is  a finite-dimensional
simple Lie superalgebra (see \cite{Zh1}). Obviously, $W (
m,n;\underline{t} ) $ is a free $ \mathcal{O} \left(
m,n;\underline{t}\right)$-module with $ \mathcal{O} \left(
m,n;\underline{t}\right) $-basis $ \{ D_r\mid r\in Y \}.$ We note
that $W ( m,n;\underline{t} ) $ possesses a \textit{standard}
$\mathbb F$-\textit{basis} $\{x^{(\alpha)}x^{u}D_r\mid \alpha\in
\mathbb{A}, u\in \mathbb{B}, r\in Y\}.$ Let $r,s\in Y$ and
$D_{rs}:\mathcal{O} ( m,n;\underline{t})\rightarrow
W(m,n;\underline{t})$ be the linear mapping  such that
\begin{equation*}
 D_{rs}(f)=(-1) ^{ \tau(r) \tau (s)}D_{r}(f)D_{s}-(-1)
^{(\tau(r)+ \tau (s)){\rm p}(f)}D_{s}(f)D_{r}\quad\mbox{for}\ f\in
\mathcal{O} ( m,n;\underline{t}).
\end{equation*}
Then the following equation holds:
\begin{equation*}
[D_{k}, D_{rs}(f)]=(-1) ^{\tau (k)\tau(r)
}D_{rs}(D_{k}(f))\quad\mbox{for}\ k,r,s \in Y; \ f\in \mathcal{O} (
m,n;\underline{t}).
\end{equation*}
Put $$ S(m,n;\underline{t}):={\rm span} _{\mathbb F}\{
D_{rs}(f)\mid r,s\in Y; \ f\in \mathcal{O}
(m,n;\underline{t})\}.$$ Then $ S(m,n;\underline{t})$ is a
finite-dimensional simple Lie superalgebra (see \cite{Zh1}). Let
$
\mathrm {div}:W(m,n;\underline{t})\rightarrow \mathcal{O}(m,n;\underline{t})$
be the \textit{divergence}  such that
$$
\mathrm {div}\big ( \sum_{r\in Y} f_{r}D_{r}\big )=\sum_{r\in Y }
(-1) ^{\tau(r) \mathrm{p}(f_{r})} D_{r}(f_{r}).
$$
Following \cite{Zh1}, put
$$
\overline{ S }(m,n;\underline{t}):=\{D\in W(m,n;\underline{t})\mid
\mathrm{div}(D)=0\}.
$$
Then $ \overline{ S }(m,n;\underline{t})$ is a subalgebra of
 $W(m,n;\underline{t})$
  and $  S (m,n;\underline{t}) $ is a subalgebra of
  $\overline{ S }(m,n;\underline{t}).$ The $\mathbb{Z}$-gradation of    $\mathcal{O}(m,n;\underline{t})$
induces naturally   a $\mathbb{Z}$-gradation structure   of
  $W(m,n;\underline{t})
=\oplus_{i=-1}^{\xi-1} W(m,n;\underline{t})_{i},$
  where
$$W(m,n;\underline{t})_{i} :={\rm span}_{\mathbb{F}}\{fD_{s}\mid
s\in Y,\ f\in \mathcal{O}(m,n;\underline{t})_{i+1}\}.$$ In
addition, $S(m,n;\underline{t})$ and
$\overline{S}(m,n;\underline{t})$ are all   $\mathbb{Z}$-graded
subalgebras of $W(m,n;\underline{t}).$ In the following sections,
$W(m,n;\underline{t}),$ $S(m,n;\underline{t}),$
$\overline{S}(m,n;\underline{t}),$ and $\mathcal{O}
(m,n;\underline{t})$ will be denoted  by $W,S, \overline{S},$ and
$\mathcal{O},$ respectively. In addition, the the even parts of
$W,$ $S$ and $\overline{S} $ will be denoted by $\mathcal{W},$
 $\mathcal{S} $ and $\overline{\mathcal{S}}, $
respectively.

 \section{Generalized Witt superalgebras}

     View $W_{\overline{1}}$ as a $\mathcal{W}$-module by
 means of the adjoint representation. In this section, our main
     purpose is to characterize the derivation space
  $\mathrm{Der}(\mathcal{W}, W_{\overline{1}}).$
    Note that the $\mathbb{Z}$-gradation of $W$ induces a $\mathbb{Z}$-gradation of
  $\mathcal{W}=\oplus_{i\geq -1}\mathcal{W}_{i}.$
 We know that gradation structures provide a powerful tool for the
study of (super)derivation algebras of Lie (super)algebra;
 in particular, the top of a $\mathbb{Z}$-graded Lie (super)algebra
 plays a predominant role (c.f. \cite{Ce,LZ1,g6}).
Following \cite{SF},
 we call
 $\mathcal{T}:=\mathrm{span}_{\mathbb{F}}\{{\Gamma}_{i}\mid i\in Y\}$ the canonical
 torus of $\mathcal{W}.$
  In the following, we first reduce every nonnegative $\mathbb{Z}$-degree
  derivation $\phi$ in $\mathrm{Der}(\mathcal{W},W_{\overline{1}})$
  to be vanishing on $\mathcal{W}_{-1};$ that is, we find an inner
  derivation $\mathrm{ad}x$ such that $(\phi-\mathrm{ad}x)(\mathcal{W}_{-1})=0, $ where
  $x\in W_{\overline{1}}.$
 In addition, we reduce the derivations in $\mathrm{Der}(\mathcal{W},W_{\overline{1}})$
  to be vanishing  on the canonical torus of
 $\mathcal{W}. $   In next step, based on these
 results,
  we shall  reduce the derivations in $\mathrm{Der}(\mathcal{W},W_{\overline{1}})$
  to be vanishing on the top $\mathcal{W}_{-1}\oplus\mathcal{W}_{0}. $

  Set
 $$\mathcal{G}:=\mathrm{span}_{\mathbb{F}}\{x^uD_i  \mid   i\in Y, u\in \mathbb{B}(n), \mathrm{p}(x^uD_i)=\overline{1}\}.$$
Note that $\mathcal{G}=C_{W_{\overline{1}}}(\mathcal{W}).$ Then $\mathcal{G}$ is a $\mathbb{Z}$-graded
  subspace of $W_{\overline{1}}.$

 In the sequel we adopt the following notation. Let $P$ be a proposition.
Define $\delta _{P}:=1$ if $P$ is true and $\delta _{P}:=0,$ otherwise. Put $\Gamma_{i}:=x_{i}D_{i}$ for $i\in Y$ and
 $\Gamma:=\sum_{i\in Y}\Gamma_{i}, $ $\Gamma':=\sum_{i\in Y_{1}}\Gamma_{i}$ and
 $\Gamma'':=\sum_{i\in Y_{0}}\Gamma_{i}.$
 We call  $\mathrm{ad}\,\Gamma$  \textit{the degree
 derivation of} $W$ (or $\mathcal{W}$), and  $\mathrm{ad}\,\Gamma'$
  and $\mathrm{ad}\,\Gamma''$ \textit{the semi-degree derivations} of $W$ (or $\mathcal{W}$).
  The following simple facts will be frequently used  in this note:
$$
\mathrm{ad}\,\Gamma (E)=rE \quad\mbox{for all}\  E\in W_{r}, \ r\in \mathbb{Z};
$$
$$
\mathrm{ad}\,\Gamma'(x^{(\alpha)}x^{u}D_{j})
=(|u|-\delta_{j\in Y_{1}})x^{(\alpha)}x^{u}D_{j}
\quad\mbox{for all}\
\alpha\in \mathbb{A}(m;\underline{t}),\, u\in \mathbb{B}(n),\, j\in Y;
$$
$$
\mathrm{ad}\,\Gamma''(x^{(\alpha)}x^{u}D_{j})
=(|\alpha|-\delta_{j\in Y_{0}})x^{(\alpha)}x^{u}D_{j}
\quad\mbox{for all}\
\alpha\in \mathbb{A}(m;\underline{t}),\, u\in \mathbb{B}(n),\, j\in Y.
$$
In particular, each standard $\mathbb{F}$-basis element $x^{(\alpha)}x^{u}D_{j}$ of $W$
is an eigenvector of the degree derivation  and the semi-degree derivation of $W.$

Similar to \cite[Lemma 2.1.1]{g6}, we have the following
 \begin{lemma}\label{lemma2.1} Suppose that $\mathcal{L}$ is a $\mathbb{Z}$-graded subalgebra of $\mathcal{W} $
  and $\mathcal{L}_{-1}=\mathcal{W}_{-1}.$
 Let  $E\in \mathcal{L}  $ and
  $\phi \in \mathrm{Der}(\mathcal{L},W_{\overline{1}}) $ satisfying
  $\phi(\mathcal{W}_{-1})=0.$
Then $\phi(E)\in \mathcal{G} $ if and only if $[E, \mathcal{W}_{-1}]\subseteq \ker\phi.$
\end{lemma}

 Analogous to \cite[Proposition 8.2, p. 192]{SF}, we have the following lemma.

  \begin{lemma} \label{lemma2.2} Let $k\leq n,\ f_{1},\ldots,f_{k}\in \Lambda(n)$
  be nonzero elements and
   ${\Gamma}_{q_{_{i}}}:=x_{q_{_{i}}}D_{q_{_{i}}},\ q_{_{i}}\in Y_{1},\ 1\leq i\leq k.$
  Suppose that

   $ \mathrm{(a)} $ \ ${\Gamma}_{q_{_{i}}}(f_j)={\Gamma}_{q_{_{j}}}(f_i) $ for
   $ 1\leq i, j \leq k;$

  $ \mathrm{(b)} $\    $ \ {\Gamma}_{q_{_{i}}}(f_i)=f_i$ for  $1\leq i\leq k.$\\
  Then there is $f\in \Lambda(n) $
  such that ${\Gamma}_{q_{_{i}}}(f)=f_i$ for $  1\leq i\leq k.$
  \end{lemma}

Analogous to \cite[Lemma 2.1.6]{g6}, we have
  \begin{lemma}\label{lemma2.3} Suppose that $\mathcal{L}$ is a $\mathbb{Z}$-graded
  subalgebra of $\mathcal{W}$ satisfying
   $\mathcal{L}_{-1}=\mathcal{W}_{-1}$. Let
  $\phi \in \mathrm{Der}(\mathcal{L},W_{\overline{1}})
  $ with
   $\mathrm{zd}(\phi)=t\geq 0.$
 Then there is $E\in (W_{\overline{1}})_{t}$
  such that
  $$\big(\phi -\mathrm{ad}E\big)( \mathcal{L}_{-1})=0.$$
  \end{lemma}

  In view of Lemma \ref{lemma2.3}, every nonnegative $\mathbb{Z}$-degree
  derivation from $\mathcal{W}$ into
  ${W}$ may be reduced to be vanishing on
  $\mathcal{W}_{-1}.$ Thus, next step is to reduce such derivations to be vanishing on the top
   $\mathcal{W}_{-1}\oplus \mathcal{W}_{0}.$  To that end,
  we first consider the canonical torus of $\mathcal{W}, $ that
  is,
   $\mathcal{T}:={\mathrm {span}}_{\mathbb{F}}\{{\Gamma}_{i}\mid i\in
   Y\}.$

 The following lemma will simplify our consideration, which tells us that
 in order to reduce derivations on the canonical torus it suffices to
  reduce these derivations on
  $\mathcal{T}':=\mathrm{span}_{\mathbb{F}}\{{\Gamma}_{j}\mid j \in Y_1\}.$

\begin{lemma}\label{lemma2.4} Suppose that
$\phi\in \mathrm{Der}_t(\mathcal{W},W_{\overline{1}}) $ with
 $ t\geq 0 $ and
 $\phi(\mathcal{W}_{-1}\cup \mathcal{T}')=0 .$
  Then $\phi({\Gamma}_{i})=0 $ for all  $ i \in Y_0.$
\end{lemma}
\begin{proof}
 (i) First consider the case $ t> 0.$ From
$\phi(\mathcal{W}_{-1})=0$ and Lemma \ref{lemma2.1} we have
 $\phi({\Gamma}_{i})\in \mathcal{G}_t $ for all $i\in Y_0.$ Thus
 one may assume that
\begin{equation}
\phi({\Gamma}_{i})=\sum_{k\in Y,\,  u\in \mathbb{B}_{t+1}}  c_{u, k}x^u
D_k \quad\mbox{where}\ c_{u, k}\in \mathbb{F}.\label{e2.1}
\end{equation}
For arbitrary $l\in Y $ and $v\in \mathbb{B}_{t+1},$ noticing that
$t+1>1, $ one may find $j\in v\backslash \{l\}.$ Clearly,
\begin{equation} [{\Gamma}_{j},\ \phi({\Gamma}_{i})]=\sum_{k\in
Y_0,\, u\in \mathbb{B}_{t+1}} c_{u, k}[{\Gamma}_{j},\ x^uD_k].
\label{e2.2}
\end{equation}
Note that
\begin{equation}[{\Gamma}_{j},\ x^uD_k]=(\delta_{j\in
u}-\delta_{jk})x^uD_k; \label{e2.3}
\end{equation}
in particular,
\begin{equation} [{\Gamma}_{j},\ x^vD_l]=x^vD_l. \label{e2.4}
\end{equation}
On the other hand,  since  $[{\Gamma}_{i},\ {\Gamma}_{j}]=0$
and $\phi({\Gamma}_{j})=0,$  one may easily see that $$[{\Gamma}_{j},\
\phi({\Gamma}_{i})]=0. $$ From (\ref{e2.2})--(\ref{e2.4})
and the equation above, we obtain that
$c_{v,l}=0.$ It follows from (\ref{e2.1}) that $\phi({\Gamma}_{i})=0 $ for all $ i\in
Y_0.$

(ii) Let us consider the case  $t=0.$
In view of Lemma \ref{lemma2.1},  we have $\phi({\Gamma}_{i})\in \mathcal{G}_0.$
 Obviously, $[\Gamma_{j},{\Gamma}_{i}]=0 $ for all $j\in Y_1,i\in Y_0.$
 Therefore,
 $$[\Gamma_{j},\phi({\Gamma}_{i})]=0 \quad\mbox{for all}\  j\in Y_1.$$
 Note that each standard basis element  of $W_{\overline{1}}$  is an  eigenvectors of
  $\mathrm{ad}\Gamma_{j}  $ for $j\in Y_1.$
  It follows from the equation displayed  above that
\begin{equation}\phi({\Gamma}_{i})=\sum_{k\in Y_1}c_{i,k}{\Gamma}_{k} \quad \mbox{where}
 \ c_{i,k}\in \mathbb{F}.\label{e2.5}
 \end{equation}
 For $j,l\in Y_1,$
 by Lemma \ref{lemma2.1}, one gets $\phi(x_jD_l)\in \mathcal{G}_0 .$
 Assume that
 $[\Gamma_{i},\phi(x_jD_l)]=\sum _{k\in Y_{1}}\lambda_{k}x_{k}D_{i}.$
 From the equation $[\Gamma_{i},x_jD_l]=0,$ we obtain that
  $$[\phi(\Gamma_{i}),x_jD_l]=-\sum _{k\in Y_{1}}\lambda_{k}x_{k}D_{i}.$$
 Then it follows from (\ref{e2.5}) that  $c_{ij}=c_{il} $ for $ j,l\in Y_1.$
 Write $c_{ij}:=c_i $  for all $j\in Y_1.$ Then
 (\ref{e2.5}) shows that
 $$\phi({\Gamma}_{i})=c_i \Gamma' \quad\mbox{for}\ i\in Y_{0}.$$
 We want to show that $ c_i=0  $ for all  $i\in Y_{0}.$
 Suppose that we are given $i\in Y_0,$ $ j,l\in Y_1.$
 Clearly, $[x_jx_lD_i,{\Gamma}_{i}]=x_jx_lD_i.$
 Applying $\phi$
 to this equation, we have
 \begin{eqnarray*}
 [\phi(x_jx_lD_i),{\Gamma}_{i}]-\phi(x_jx_lD_i)&=&-[x_jx_lD_i,\phi({\Gamma}_{i})]\\
 &=&-[x_jx_lD_i, c_{i}{\Gamma}']\\
 &=&2c_{i} x_jx_lD_i.
 \end{eqnarray*}
 By Lemma \ref{lemma2.1},
 it is easily seen that $\phi(x_jx_lD_i)\in \mathcal{G}_1. $ Thus
 one may assume that
  $\phi(x_jx_lD_i)=\sum_{k\in Y,\,    u\in \mathbb{B}_{2}}c_{u,k}x^uD_k.$
 Note that $[x^uD_k,x_iD_i]=\delta_{ki}x^uD_i.$  It follows that
 $$\sum_{k\in Y,\, u\in \mathbb{B}_{2}}(\delta_{ki}-1)c_{u,k}x^uD_k=2c_{i} x_jx_lD_i. $$
 A comparison of the coefficients of  $x_jx_lD_i$ in the equation above yields that
  $2c_{i}=0 $ for $i\in Y_0.$
Since $\mathrm{char  }\mathbb{F}\neq 2,$ we have $c_{i} =0 $ for
all $i\in Y_0.$ So far, we have proved that $\phi({\Gamma}_{i})=0
$ for all $ i\in Y_0.$

 Now, by (i) and (ii), we obtain the desired result.
 \end{proof}

We first consider the homogeneous derivations of odd positive
$\mathbb{Z}$-degree.

\begin{lemma} \label{lemma2.5} Suppose that
$\phi\in \mathrm{Der}_{t}(\mathcal{W},W_{\overline{1}}) $ where
$\mathrm{zd}(\phi)=t\geq 1$ is odd. If $\phi(\mathcal{W}_{-1})=0,
$ then there is $z\in \mathcal{G}_t   $ such that
$(\phi-\mathrm{ad}z)(\mathcal{T}')=0. $
\end{lemma}
\begin{proof}  Using Lemma \ref{lemma2.1} and noting that $t$ is odd,  we
may assume that
\begin{equation}
\phi(\Gamma_{i})=\sum_{r\in Y_1}f_{ri}D_r
\quad \mbox{where}\ i\in Y_1,\ f_{ri}\in\Lambda (n).
\label{e2.14}
\end{equation}
Applying $\phi$
to the equation that $[\Gamma_{i}, \Gamma_{j}]=0$ for $ i, j\in Y_1, $
we have
$$\sum_{r\in Y_1}(\Gamma_{i}(f_{rj})-\Gamma_{j}(f_{ri}))D_r+f_{ji}D_j-f_{ij}D_i=0.$$
Consequently,
\begin{equation}\Gamma_{i}(f_{rj})=\Gamma_{j}(f_{ri})
 \quad \mbox{whenever}\ r\neq i, j;\label{e2.15}
\end{equation}
\begin{equation}\Gamma_{j}(f_{ii})=\Gamma_{i}(f_{ij})-f_{ij}
\quad \mbox{whenever}\ i\neq j.\label{e2.16}
\end{equation}
For $r, i\in Y_1,$
one may assume that $f_{ri}=\sum_{|u|=t+1}c_{uri}x^u, \  c_{uri}\in \mathbb{F}.$
By (\ref{e2.15}), we have
$$c_{uri} \delta _{j\in u}=c_{urj} \delta _{i\in u}
 \quad \mbox{whenever}\ r\neq i, j.$$
This implies that
$$c_{uri} \neq 0\ \mbox{and }\ j\in u \Longleftrightarrow
c_{urj} \neq 0\ \mbox{and }\ i\in u.$$
Let $r\neq i $ and assume that
 $c_{uri} \neq 0.$
Then the implication relation above shows that $i\in u.$
Accordingly,
\begin{equation}
\Gamma_{i}(f_{ri})=f_{ri} \quad \mbox{whenever}\ r\neq i.\label{e2.17}
\end{equation}
For any fixed $r\in Y_1, $
Lemma \ref{lemma2.2} ensures that there is $\overline{f}_r\in \Lambda(n)$
such that
\begin{equation}
\Gamma_{i}(\overline{f}_{r})=f_{ri}
\quad \mbox{for all}\  i\in Y_1\backslash \{r\}.\label{e2.18}
\end{equation}
Assert that
\begin{equation}
\Gamma_{i}(f_{ii})=0  \quad \mbox{for all}\  i\in Y_1.\label{e2.19}
\end{equation}
Using (\ref{e2.16}) and noticing the fact that ${\Gamma}^2_i=\Gamma_{i},$
we obtain that
$$\Gamma_{j}\Gamma_{i}(f_{ii})=\Gamma_{i}\Gamma_{j}(f_{ii})
={\Gamma}^2_i(f_{ij})-\Gamma_{i}(f_{ij})=0  \quad\mbox{for}\ j\neq i.$$
Note that $\mathrm{zd}(f_{ii})=t+1\geq 2  $
and $\Gamma_{i}(x^u)=\delta_{i\in u}x^u.$  (\ref{e2.19}) follows.

 For $r\in Y_1,$
 put $f_r:=-f_{rr}+{\Gamma}_r(\overline{f}_r).$
 Obviously, $f_r\in \Lambda(n).  $
 It follows from (\ref{e2.19}) that
 \begin{equation}
 {\Gamma}_r(f_r)-f_r=-{\Gamma}_r(f_{rr})+{\Gamma}_r^{2}(\overline{f}_r)+f_{rr}-
 {\Gamma}_r(\overline{f}_r)=f_{rr}.
 \label{e2.20}
 \end{equation}
 For $i\in Y_1 \backslash r, $
  by (\ref{e2.16}) and (\ref{e2.18}) we obtain that
 \begin{eqnarray*}
 \Gamma_{i}(f_r)&=&-\Gamma_{i}(f_{rr})+\Gamma_{i}{\Gamma}_r(\overline{f}_r)
 =-\Gamma_{i}(f_{rr})+{\Gamma}_r\Gamma_{i}(\overline{f}_r)\\
 &=&-({\Gamma}_r(f_{ri})-f_{ri})+{\Gamma}_r(f_{ri})
 =f_{ri}.
 \end{eqnarray*}
 Let $z':=-\sum_{r\in Y_1}f_rD_r.$
A combination of  (\ref{e2.20}) and the equation above yields that
 for $i\in Y_1,$
\begin{eqnarray*}
 [z', \Gamma_{i}]&=&-\sum_{r\in Y_1}[f_rD_r, \Gamma_{i}]
 =\sum_{r\in Y_1}\Gamma_{i}(f_r)D_r-f_iD_i\\
 &=&\sum_{r\in Y_1\backslash i}\Gamma_{i}(f_r)D_r+(\Gamma_{i}(f_i)-f_i)D_i\\
 &=&\sum_{r\in Y_1\backslash i}f_{ri}D_r+f_{ii}D_i
 =\phi(\Gamma_{i}).
 \end{eqnarray*}
Let $z$ be the $t$-component of $z'.$ Since $\mathrm{zd}(\phi)=t,$ one gets
$[z, \Gamma_{i}]=\phi(\Gamma_{i}) $ for all $i\in Y_1.$
Putting $\psi:=\phi-\mathrm{ad}z, $ then $\psi\in \mathrm{Der}_t
(\mathcal{W},W_{\overline{1}}) $ and $\psi(\Gamma_{i})=0 $ for all $i\in Y_1.$
\end{proof}

  \begin{lemma}\label{lemma2.6}
  Suppose that $\phi\in \mathrm{Der}(\mathcal{W},W_{\overline{1}}) $ and $ \mathrm{zd}(\phi)=t\geq 0 $
is even. If $\phi(\mathcal{W}_{-1})=0 $ then there is $z\in
\mathcal{G}_t$ such that $(\phi-\mathrm{ad}z)( \mathcal{T}')=0. $
\end{lemma}

\begin{proof}
 Since $\mathrm{zd}(\phi)=t$
is even, by Lemma \ref{lemma2.1}, one may assume that
\begin{equation}
\phi(\Gamma_{i})=\sum_{r\in Y_0}f_{ri}D_r
\quad \mbox{where}\ i\in Y_1,\ f_{ri}\in\Lambda (n).
\label{e2.32}
\end{equation}
Analogous to the proof of Lemma \ref{lemma2.5},  one may easily
show that
\begin{equation}\Gamma_{i}(f_{rj})=\Gamma_{j}(f_{ri})
 \quad \mbox{for all}\ i, j\in Y_1.\label{e2.33}
\end{equation}
Suppose that
\begin{equation}
f_{ri}=\sum_{u\in \mathbb{B}_{t+1}}c_{u, r, i}x^u
\quad\mbox{where}\ c_{u, r, i}\in \mathbb{F}.\label{e2.34}
\end{equation}
Then we obtain from  (\ref{e2.33}) and (\ref{e2.34}) that
$$c_{u, r, i}\delta _{j\in u}=c_{u, r, j}\delta _{i\in u}  \quad\mbox{for all}\
 i, j\in Y_1,\
 r\in Y_0. $$
Consequently, for $i, j\in Y_1,$ $ r\in Y_0 $ and $ u\in \mathbb{B}_{t+1},$
\begin{equation}
c_{u, r, i}\neq 0 \ \mbox{and}\ j\in u\Longleftrightarrow
c_{u, r, j}\neq 0 \ \mbox{and}\ i\in u.\label{e2.35}
\end{equation}

Let us complete the proof of this lemma.
Assume that $\mathrm{zd}(\phi)=t\geq 2.$
If $c_{u, r, i}\neq 0 $ for  $i\in Y_1,$
one may pick $j\in u\backslash i.$
By (\ref{e2.35}), we have  $i\in u.$ Assume that $\mathrm{zd}(\phi)=0 .$ Then
(\ref{e2.35}) implies that there is at most one nonzero summand
$c_{\langle i\rangle, r,i}x_{i} $ in the right-hand side of  (\ref{e2.34}).
Summarizing,   every nonzero summand in the right-hand
side of (\ref{e2.34}) possesses the factor $x_i.$
Therefore,
\begin{equation}
\Gamma_{i}(f_{ri})=f_{ri}  \quad\mbox{for all}\ \ i\in Y_1,\ r\in Y_0.\label{e2.36}
\end{equation}
For any fixed $r\in Y_0,$
by (\ref{e2.33}) and (\ref{e2.36}), $\{f_{r, m+1 }, f_{r,m+2}, \ldots, f_{r,m+n}\}$
fulfills the conditions of Lemma \ref{lemma2.2}.  Hence, there is
$f_r\in \Lambda(n) $
such that
$$\Gamma_{i} (f_{r})=f_{ri}\quad \mbox{for}\ i\in Y_1.$$
Let $z':=-\sum_{r\in Y_0}f_rD_r.$
Then (\ref{e2.32}) and the equation above show  that
  $[z', \Gamma_{i}]=  \phi(\Gamma_{i}) $ for $i\in Y_1. $
Let $z$ be $t$-component of $z'.$
Then $z\in \mathcal{G}_t $ and
$(\phi-\mathrm{ad}z)(\Gamma_{i})=0 $ for all $i\in Y_1.$
\end{proof}

 Now we come to  the following main result.

 \begin{proposition} \label{proposition2.7}
 Let $\phi$ be a homogeneous derivation from $\mathcal{W}$ into $W_{\overline{1}}$ with
 nonnegative $\mathbb{Z}$-degree $t.$ Then
  $\phi$ can be reduced to be vanishing on $\mathcal{W}_{-1}$ and the canonical
  torus of $\mathcal{W};$
 that is, there is $E\in (W_{\overline{1}})_{t} $ such that
 $(\phi-\mathrm{ad}E)\mid _{_{\mathcal{W}_{_{-1}}+\mathcal{T}}}=0.$
 \end{proposition}

 \begin{proof} By Lemma \ref{lemma2.3}, there is $E'\in (W_{\overline{1}})_{t} $
 such that $(\phi-\mathrm{ad}E')(\mathcal{W}_{-1})=0.$
 Then by Lemmas \ref{lemma2.4}---\ref{lemma2.6}, there is
 $E''\in \mathcal{G}_{t} $ such that
 $(\phi-\mathrm{ad}E'-\mathrm{ad}E'')(\mathcal{T})=0.$
 Putting $E:=E'+E'',$
 then $(\phi-\mathrm{ad}E)(\mathcal{W}_{-1}+\mathcal{T})=0.$\end{proof}

  In the following, using Proposition \ref{proposition2.7}, we  first  reduce
 every nonnegative $\mathbb{Z}$-degree derivation from $\mathcal{W}$ into $W_{\overline{1}}$ to be vanishing on the top
 $\mathcal{W}_{-1}\oplus\mathcal{W}_{0}$ of $\mathcal{W};$ then we determine
 the $\mathbb{Z}$-homogeneous components
 $\mathrm{Der}_{t}(\mathcal{W},W_{\overline{1}}) $ for  $t\geq 0. $

 \begin{proposition}\label{proposition3.1}
 Let $\phi\in \mathrm{Der}_{t}(\mathcal{W},W_{\overline{1}}) $ with $t\geq 0.$ Then there is
 $E\in (W_{\overline{1}})_{t} $ such that
 $(\phi-\mathrm{ad}E)\mid _{_{\mathcal{W}_{_{-1}}\oplus\mathcal{W}_{_{0}}}}=0.$
 \end{proposition}

 \begin{proof} By Proposition \ref{proposition2.7}, without loss of generality
 we may assume that $\phi(\mathcal{W}_{-1}+\mathcal{T})=0.$

 (i)  We first consider $\phi(x_kD_i)  $ where $i,k\in Y_0$
 with $i\neq k.$ By Lemma \ref{lemma2.1},
 $\phi(x_kD_i)\in \mathcal{G}_t.$ Assume that
 $\phi(x_kD_i)=\sum_{r\in Y, u\in \mathbb{B}_{t+1} } c_{u,r}x^{u}D_{r}$
 where $c_{u,r}\in \mathbb{F}.$
  If $t$ is even, then
  $\phi(x_kD_i)=\sum_{r\in Y_{0}, u\in \mathbb{B}_{t+1} } c_{u,r}x^{u}D_{r}.$
  Note that
 $[\Gamma_{j}, \phi(x_kD_i)]=0$ for arbitrary $j\in Y_{1}.$ It
 follows that $\phi(x_kD_i)=0.$ If $t$ is odd, then
 $\phi(x_kD_i)=\sum_{r\in Y_{1}, u\in \mathbb{B}_{t+1} } c_{u,r}x^{u}D_{r}.$
Then
\begin{equation*}
\phi(x_kD_i)=\phi([\Gamma_{k},x_kD_i])=[\Gamma_{k},\phi(x_kD_i)]=0.
\end{equation*}

 (ii) We next consider  $\phi(x_kD_l)$ where $k,l\in Y_1$ with $k\neq l.$

 (a) Suppose that $t$ is even. Just as in (i) one may assume that
$$\phi(x_kD_l)=
\sum_{r\in Y_{0}, u\in \mathbb{B}_{t+1} } c_{u,r}x^{u}D_{r}\quad \mbox{where}\
c_{u,r}\in \mathbb{F}.
$$ Then,
from the equation that $[\Gamma_{i},\phi(x_kD_l)]=0$ for all $i\in Y_{0}$,
one  gets $c_{u,r}=0$ for all $r\in Y_{0},$ $ u\in \mathbb{B}_{t+1}.$ Hence,
$\phi(x_kD_l)=0.$

 (b) Suppose that  $t $ is odd.  We proceed in two cases
 $t\geq 3$ and $1\leq t\leq 2$
  to show that $\phi(x_kD_l)=0 $
 for $k,l\in Y_1.$
 Suppose that  $t\geq 3.$ By Lemma \ref{lemma2.1}, one may assume that
 \begin{equation}\phi(x_kD_l)=\sum_{r\in Y_{1},\ \mid u\mid \geq 4} c_{u,r}x^{u}D_{r}
 \quad \mbox{where}\  c_{u,r}\in \mathbb{F}.\label{e3.2}
 \end{equation}
 Given $v\in \mathbb{B} $ with $\mid v\mid \geq 4 $  and  $s\in Y_{1},$  choose
 $q\in v\setminus \{k,l,s\}.$  Then $[\Gamma_{q}, x^{v}D_{s}]=x^{v}D_{s}.$
 On the other hand, since  $[\Gamma_{q},x_kD_l]=0,$ we have
 $[\Gamma_{q},\phi(x_kD_l)]=0.$ Note that each standard basis
 element of $W$ is an eigenvector of
 $\Gamma_{q}$ and
 $[\Gamma_{q}, x^{v}D_{s}]=x^{v}D_{s}.$  It follows from (\ref{e3.2}) that
 $c_{v,s}=0 .$ Therefore, $\phi(x_kD_l)=0.$

 Finally we consider the case $t=1.$
 Clearly,
 $[\Gamma',x_kD_l]=0.$
Consequently,
  $[\Gamma',\phi(x_kD_l)]=0.$
On the other hand,  by Lemma \ref{lemma2.1},  $\phi(x_kD_l)\in \mathcal{G}. $ Thus
\begin{equation*}
0=[\Gamma',\phi(x_kD_l)]=\phi(x_kD_l).
\end{equation*}
The proof is complete.
\end{proof}

 In order to determine the homogeneous derivation subspace
 $\mathrm{Der}_{t}(\mathcal{W},W_{\overline{1}}) $ for  $t\geq 0,$
 we need the  generator set of $\mathcal{W} $ (see \cite[Proposition 2.2.1]{g6}).

 \begin{lemma} \label{lemma3.2} $\mathcal{W}$ is generated by
 $\mathcal{P}\cup \mathcal{N}\cup \mathcal{M},$
  where
  $$\mathcal{P}:=\{x_k x_lD_{i}\mid k,l\in Y_1,\ i\in Y_0\},$$
  $$\mathcal{N}:=\{x_k x_lD_{i}\mid  k\in Y_0,\ l,i\in Y_1\},$$
  $$\mathcal{M}:=\{x^{(k\varepsilon_{i})}D_{j}\mid 0\leq k\leq \pi_{i},\ i,j\in Y_0\}.$$
  \end{lemma}

Now we can determine the homogeneous derivations of
nonnegative $\mathbb{Z}$-degree:

 \begin{proposition}\label{proposition3.3} Let $t\in \mathbb{N}_{0}.$
 Then $\mathrm{Der}_t (\mathcal{W}, W_{\overline{1}})
 =\mathrm{ad} (W_{\overline{1}})_t.$
\end{proposition}
 \begin{proof}  If suffices to show the inclusion ``$\subset$".
 Let $\phi\in \mathrm{Der}_t (\mathcal{W}, W)$. By Proposition \ref{proposition3.1},
 one may assume that $\phi(\mathcal{W}_{-1}\oplus\mathcal{W}_{0})=0. $
 In the following we consider the application of  $\phi$ to
 $\mathcal{P},$ $\mathcal{N}$ and  $\mathcal{M},$ respectively.

 (i) First consider $\mathcal{P}.$ Let $i,j\in Y_1,$ $k\in Y_0.$ In view of Lemma
 \ref{lemma2.1},  we have
 $\phi(x_{i}x_jD_k)\in \mathcal{G}_{t+1}. $
 Clearly, $[\Gamma'',x_{i}x_jD_k]=-x_{i}x_jD_k.$
 Applying $\phi,$
 we obtain that
 $$[\Gamma'' ,\phi(x_{i}x_jD_k)]=-\phi (x_{i}x_jD_k). $$
 It follows that $\phi (x_{i}x_jD_k)$ is of the form:
 \begin{equation}\phi (x_{i}x_jD_k)
 =\sum_{r\in Y_{0},\ u\in \mathbb{B}_{t+2}}c_{u,r}x^{u}D_{r}
 \quad\mbox{where}\ c_{u,r}\in \mathbb{F}.\label{e3.3}
 \end{equation}
 Since $n>2,$
picking  $ l\in Y_1 \setminus \{i,j\},$
 one gets $[x_ix_jD_k,x_lD_l]=0.$
 It follows that $[\phi(x_ix_jD_k),x_lD_l]=0 .$
 Furthermore,
 \begin{equation}
 \Big[\sum_{r\in Y_0,\ u\in \mathbb{B}_{t+2}}
 c _{u,r}x^uD_r,x_lD_l\Big]=0  \quad \mbox{for all}\    l\neq i,j.\label{e3.4}
 \end{equation}
 Now it follows from (\ref{e3.3}) and (\ref{e3.4}) that $c _{u,r}=0 $
 unless $u=\{i,j\}.$
 Thus  (\ref{e3.3}) gives
 \begin{equation} \phi(x_ix_jD_k)
 =\sum_{r\in Y_0,u\in \mathbb{B}\atop u =\{i,j\}}c_{u,r}x^{u}D_r.\label{e3.5}
 \end{equation}
 If $t>0 $ then (\ref{e3.5}) implies that  $\phi(x_ix_jD_k)=0;$
if $t=0 $, we also obtain from (\ref{e3.5}) that
$\phi(x_ix_jD_k)=0, $  since $\phi\in \mathrm{Der}(\mathcal{W},
W_{\overline{1}}).$

(ii) Let us show that $\phi(\mathcal{N})=0.$ Let $i\in Y_{0},$ $j,
k\in Y_{1}.$ By Lemma \ref{e2.1}, $\phi(x_{i}x_{j}D_{k})\in
\mathcal{G}_{t+1}.$ Since $[\Gamma'',
x_{i}x_{j}D_{k}]=x_{i}x_{j}D_{k},$ as in (i) one may assume that
$$\phi(x_{i}x_{j}D_{k}) =\sum_{r\in Y_{0}, u\in \mathbb{B}_{t+2}}c_{u,
r}x^{u}D_{r}
 \quad\mbox{where}\ c_{u, r}\in \mathbb{F}.$$ Then
$$-\phi(x_{i}x_{j}D_{k})=[\Gamma'', \phi(x_{i}x_{j}D_{k})]=\phi(x_{i}x_{j}D_{k}).$$
Since $\mathrm{char}\mathbb{F}\neq 2,$ it follows that
 $$\phi(x_{i}x_{j}D_{k})=0 \quad\mbox{for all}\ i\in Y_{0},\ j, k\in Y_{1};$$
that is, $\phi(\mathcal{N})=0.$

 (iii) Just as in the proof of \cite[Lemma 3.1.4]{g6}, one may show
 that
   $\phi(\mathcal{M})=0.$

Now, Lemma \ref{lemma3.2} shows that $\phi=0.$
\end{proof}

 In view of Proposition \ref{proposition3.3},  in order to determine the
 derivation space $\mathrm{Der}(\mathcal{W}, W_{\overline{1}})$
 it suffices to determine the derivations of  negative $\mathbb{Z}$-degree. We first
 consider the derivations of $\mathbb{Z}$-degree $-1.$

 \begin{lemma}\label{lemma4.1}  Suppose that
 $\varphi\in \mathrm{Der}_{-1}(\mathcal{W},W_{\overline{1}}) $ and
  $\varphi(\mathcal{W}_0)=0.$
 Then  $\varphi=0.$
 \end{lemma}
\begin{proof} We first assert that $\varphi(\mathcal{N})=\varphi(\mathcal{P})=0.$
 Given $i,j\in Y_1,k\in Y_0, $ by Lemma \ref{lemma2.1},
 $\varphi (x_ix_jD_k)\in \mathcal{G}_0.$ Thus one may assume that
 $\varphi (x_ix_jD_k)=\sum_{r\in Y_{1},\ s\in Y_{0}}c_{r,s}x_{r}D_{s},$
 $c_{r,s}\in \mathbb{F}.$
  Then,
 since $[\Gamma', x_ix_jD_k]=2x_ix_jD_k,$
 we have
  $$\varphi(x_ix_jD_k)=[\Gamma', \varphi(x_ix_jD_k)]=2\varphi(x_ix_jD_k).$$
   It follows that $\varphi(x_ix_jD_k)=0;$ that is, $\varphi(\mathcal{P})=0.$
  Similarly,  applying $\varphi$
 to the equation $[\Gamma'',x_kx_{i}D_j]=x_kx_iD_j $ gives
 $$-\varphi( x_kx_{i}D_j)=[\Gamma'', \varphi( x_kx_{i}D_j)]=\varphi( x_kx_{i}D_j).$$
 It follows that $\varphi( x_kx_{i}D_j)=0. $
 Hence,  $\varphi(\mathcal{N})=0.$

 It remains to show that $\varphi(\mathcal{M})=0.$ Given $k\in Y_0,$
just as in the proof of \cite[Lemma 3.2.6]{g6}, one may prove by
induction on $r$ that
 \begin{equation*}\varphi(x^{(r\varepsilon_k)}D_k)
 =0 \quad \mbox{for all}\ r\in \mathbb{N}.
 \end{equation*}
  From this one may easily prove that
   $\varphi(\mathcal{M})=0.$
 Summarizing, by Lemma \ref{lemma3.2}, $\varphi=0.$
\end{proof}

Now we can determine the derivations  from $\mathcal{W}$ into $W_{\overline{1}} $ of $ \mathbb{Z}$-degree $-1.$

 \begin{proposition}\label{proposition4.2}
  $ \mathrm{Der}_{-1}(\mathcal{W}, W_{\overline{1}})
  =\mathrm{ad} (W_{\overline{1}})_{-1}.$
 \end{proposition}

 \begin{proof}  The inclusion ``$\supset$'' is clear.
 Let $\phi \in \mathrm{Der}_{-1}(\mathcal{W}, W).$
 For $i\in Y_0,k\in Y_1,$
 applying $\phi$
 to the equation that $[\Gamma_{i}, \Gamma_{k}]=0,$
  we have
$
[\phi(\Gamma_{i}), \Gamma_{k}]+[ \Gamma_{i},\phi(\Gamma_{k}
)]=0.
$
  As $\phi(\Gamma_{i}),\ \phi(\Gamma_{k})\in  {W}_{-1}\cap W_{\overline{1}},$
  we have $[ \Gamma_{i},\phi(\Gamma_{k} )]=0 $ and therefore,
  $[\phi(\Gamma_{i}), \Gamma_{k}]=0 $ for all $k\in Y_1.$ This
  implies that $\phi(\Gamma_{i})=0$ for $i\in Y_0.$  It follows that
  $\phi(x_{i}D_{j})=0$ for all $i,j\in Y_{0}.$

 For $k\in Y_1,$ just as in the proof of \cite[Proposition
 3.2.7]{g6}, one may prove that there are $c_{k}\in \mathbb{F} $ such that
  $\phi(\Gamma_{k})=c_{k}D_k$  and
 $
 \phi(x_kD_l)=c_{k}D_l$ for all $k,l\in Y_{1}.
 $
By Lemma \ref{lemma4.1},
$\phi=\sum_{r\in Y_{1}} c_{r}D_r \in \mathrm{ad}
 (W_{\overline{1}})_{-1}.$
\end{proof}

Analogous to \cite[Lemma 3.2.8]{g6}, we also have the following

 \begin{lemma} \label{lemma4.3}
  Let $\phi \in \mathrm{Der}_{-q}(\mathcal{W}, W_{\overline{1}}) $
 with $q>1.$
 If $\phi (x^{(q\varepsilon_i)}D_i)=0 $ for all $i\in Y_0,$
  then $\phi=0.$
  \end{lemma}

\begin{proposition}\label{proposition4.4}  Suppose that $q>1. $
Then $\mathrm{Der}_{-q}(\mathcal{W},W_{\overline{1}})=0.$
\end{proposition}
 \begin{proof} Let $ \phi\in
\mathrm{Der}_{-t}(\mathcal{W},W_{\bar{1}}).$ In view of Lemma
\ref{lemma4.3}, it is sufficient to show that
 $$ \phi(x^{(q\varepsilon_i)}D_{i})=0 \quad \mbox {for
 all}\; i\in Y_0.  $$
Note that $[
 \Gamma^{'}, x^{ (q\varepsilon_i)}D_{i} ]=0 $ and
 $\phi(x^{(q\varepsilon_i)}D_{i})\in
(W_{\bar{1}})_{-1}.$ It follows that
 $$ 0=[ \Gamma^{'}, \phi(x^{(q\varepsilon_i)}D_{i})]
 =-\phi( x^{(q\varepsilon_i)}D_{i}).$$
 The proof is complete.
\end{proof}

 \begin{theorem}\label{theorem4.6}  $\mathrm{Der}(\mathcal{W}, W_{\overline{1}})=\mathrm{ad}(W_{\overline{1}}).$
  \end{theorem}
 \begin{proof}   By Propositions
  \ref{proposition3.3}, \ref{proposition4.2} and \ref{proposition4.4},
  ``$\subset$"  holds.  The converse inclusion is clear.
  \end{proof}

   By \cite[Theorem 2]{ZZ1}, \cite[Theorem 3.2.11]{g6} and Theorem \ref{theorem4.6},
   the even part and the odd part of the superderivation algebra
   of the finite-dimensional generalized Witt superalgebra $W$ coincide with
    the derivation algebra of the even part of $W$ and
    the derivation space of the even part into the odd part of  $W, $ respectively;
 that is,
    $(\mathrm{Der}\, W )_{\overline{0}}
   =\mathrm{Der}(W_{\overline{0}}),$
   $(\mathrm{Der}\, W )_{\overline{1}}
   =\mathrm{Der}(W_{\overline{0}},\, W_{\overline{1}}). $

\section{Special superalgebras}

Recall the canonical torus $ \mathcal{T}_{\mathcal{S}} $ of $ \mathcal{S} $ (c.f. \cite{g6}). Clearly,
 $$\{ x_rD_r-x_sD_s | \tau(r)=\tau(s); r,s\in Y\} \cup \{x_rD_r+x_sD_s| \tau(r)\neq \tau(s); r,s\in
 Y\}$$
 is an  $ \mathbb{F}$-basis   of $  \mathcal{T}_{\mathcal{S}} $ consisting of  toral
 elements.

 The following fact is simple.
 \begin{lemma}\label{gt2.1.7}

$\mathcal{S}_0=\mathrm{span}_{\mathbb{F}}\{\mathcal{T}_{\mathcal{S}}\cup\{x_rD_s
\,|\, \tau(r)=\tau(s),r\neq s; r,s\in Y \}\}. $
 \end{lemma}

 Put
$$\mathcal{Q}:=\{D_{ij}(x^{(r\varepsilon_j)})\, |\, i,j\in Y_0,r\in
 \mathbb{N}_0\};$$
 $$\mathcal{R}:=\{ D_{il}(x^{(2\varepsilon_i)}x_k) \,|\, i\in Y_0,k,l\in Y_1 \}\cup\{ D_{ij}(x_ix^{v})\,
 |\, i,j\in Y_0,v\in \mathbb{B}_2\}.$$

 We need the generator set of $  \mathcal{S} $ (see \cite[Proposition 2.2.3]{g6}).

\begin{lemma}\label{gt2.1.8}

 $  \mathcal{S} $ is generated by
 $ \mathcal{Q}\cup\mathcal{R}\cup\mathcal{S}_0. $
\end{lemma}

  In the following we consider the top of $\mathcal{S}.$

\begin{lemma}\label{gt3.1.1}

Suppose that $  \phi\in \mathrm{Der}(\mathcal{S},W_{\bar{1}})$
with $ \mathrm{zd}(\phi)\geq 0$ and that
 $  \phi(\mathcal{S}_{-1}+\mathcal{S}_0)=0 .$
 Then
 \begin{itemize}
\item[$\mathrm{(i)}$] $\phi(D_{il}(x^{(2\varepsilon_i)}x_k))=0 \;
\mbox {for all}\; i\in Y_0,\, k,l\in Y_1.
$
\item[$\mathrm{(ii)}$] $  \phi(D_{ij}(x_ix^{v}))=0 $  for all
 $ i,j\in Y_0 $ and   $
v\in \mathbb{B}_2.$
\item[$\mathrm{(iii)}$]$\phi(D_{ij}(x^{(a\varepsilon_i)}))=0\; $ for
all $i,j\in Y_0 $ and all $a\in
\mathbb{N}.$
\end{itemize}
 \end{lemma}
\begin{proof} (i) The proof is similar to the one of
\cite[Lemma 4.1.1]{g6}. Our discussion here for $\mathrm{zd}(\phi)$ being odd is
completely analogous  to
one in \cite[Lemma 4.1.1]{g6} for  $\mathrm{zd}(\phi)$ being even; and, the
discussion here for $\mathrm{zd}(\phi)$ being even is completely analogous to
one in \cite[Lemma 4.1.1]{g6} for $\mathrm{zd}(\phi)$ being odd.

Similar to \cite[Lemmas 4.1.2, 4.1.3]{g6}, one may prove (ii) and (iii) in the same way.
\end{proof}

As a direct consequence of Lemmas \ref{gt2.1.8} and \ref{gt3.1.1}, we have the following.

\begin{corollary} \label{gt3.1.4}

 Suppose that $ \phi \in \mathrm{Der}(\mathcal{S}, W_{\bar{1}}) $  with
$ \mathrm{z}d(\phi)\geq 0 $ and that $
\phi(\mathcal{S}_{-1}+\mathcal{S}_0)=0.$  Then $ \phi=0$.
\end{corollary}

In order to describe  the derivations of nonnegative degree we first give two technical lemmas which will simplify
our discussion.

\begin{lemma}\label{gt3.2.1} 

Suppose that $ \phi\in \mathrm{Der}_t(\mathcal{S},W_{\bar{1}})$
and
$\phi(\mathcal{S}_{-1})=0.$ \\
 $\mathrm{(i)}$ If $ t=n-1 $ is even, then $ \phi(\Gamma_{1^{'}}-\Gamma_k)=0 $
 \,for all \,$  k\in Y_1\setminus 1^{'}.$ \\
 $\mathrm{(ii)}$ If $ t=n-1 $ is odd, then there is  $  \lambda\in \mathbb{F} $ such that
$$(\phi-\lambda
\mathrm{ad}(x^{\omega}D_{1^{'}}))(\Gamma_{1^{'}}-\Gamma_k)=0
 \; \mbox {for all} \; k\in Y_1\setminus 1^{'}.$$
$\mathrm{(iii)}$ If $ t > n-1,$ then $ \phi=0.$
\end{lemma}
\begin{proof} (i) The proof is completely analogous to the
one of  \cite[Lemmas 4.2.1(i)]{g6}.

    (ii) The proof is completely analogous to the
one of  \cite[Lemmas 4.2.1(ii)]{g6}.

   (iii)
   Using Lemma \ref{lemma2.1} and   induction on $r$ one may easily prove that
   $\phi(S_r)=0$ for all $r\in \mathbb{N}.$
   \end{proof}

Analogous to \cite[Lemmas 4.2.2]{g6}, we have

\begin{lemma}\label{gt3.2.2} 

Suppose that $ \phi\in \mathrm{Der}(\mathcal{S},W_{\bar{1}})$ and
$
\mathrm{zd}(\phi)\geq 0 $ is even. \\
   $\mathrm{(i)}$ If $\mathrm{zd}(\phi)<n-1$ and
   $$\phi(\Gamma_{1'}-\Gamma_{2'})=\phi(\Gamma_{1'}-\Gamma_{3'})=\cdots
   =\phi(\Gamma_{1'}-\Gamma_{n'})=0,$$ then
   $$\phi(\Gamma_{1}-\Gamma_{2})=\phi(\Gamma_{1}-\Gamma_{3})=\cdots
   =\phi(\Gamma_{1}-\Gamma_{m})=0;\quad \phi(\Gamma_{1}+\Gamma_{1'})=0. $$

  $\mathrm{(ii)}$ If $\mathrm{zd}(\phi)=n-1,$ then there are
   $\lambda_1,\ldots,\lambda_m\in \mathbb{ F}  $ such that
   $$\Big(\phi-\mathrm{ad}\Big(\sum _{i\in Y_{0}}\lambda_ix^{\omega} D_i\Big)\Big)
   (\Gamma_1-\Gamma_j)=0
   \quad \mbox{for all} \; j\in Y_0\setminus 1; $$
   $$\Big(\phi-\mathrm{ad}\Big(\sum _{i\in Y_0}\lambda_ix^{\omega} D_i\Big)\Big)
   (\Gamma_1+\Gamma_{1'})=0.
   \eqno$$
\end{lemma}

 Recall the canonical torus of $ \mathcal{S} $
 $$ \mathcal{T}_{\mathcal{S}} := \mathrm{span}_\mathbb{F}
 \{ \Gamma_1-\Gamma_2,\ldots,\Gamma_1-\Gamma_m,\ldots,\Gamma_1+\Gamma_{1^{'}},
 \Gamma_{1^{'}}-\Gamma_{2^{'}},\ldots,\Gamma_{1^{'}}-\Gamma_{n^{'}} \}. $$

As a direct consequence of Lemma \ref{gt3.2.1} (iii)
    and Lemma \ref{gt3.2.2}, we have the following fact:

\begin{corollary}\label{gt3.2.3}

Suppose that
   $\phi\in \mathrm{Der}(\mathcal{S},W_{\bar{1}})$ is homogeneous derivation
   of nonnegative even $\mathbb{Z}$-degree such that $\phi(\mathcal{S}_{-1})=0$ and
   $\phi(\Gamma_{k}-\Gamma_{1'})=0$ for all $k\in Y_1\setminus 1'.$ Then there is
   $E\in \mathcal{G}$
   such that
   $(\phi-\mathrm{ad}E) $ vanishes on the canonical torus $\mathcal T_{\mathcal S}.$
\end{corollary}

 Now we  prove the following two key lemmas. First, consider
 the derivations of even $\mathbb{Z}$-degree.
\begin{lemma}\label{gt3.2.4}

Suppose that
   $\phi \in \mathrm{Der}_{t}(\mathcal{S},W_{\bar{1}}),$
   where $t\geq 0$ is even.  If $\phi(\mathcal{S}_{-1})=0,$
   then there is $D\in \mathcal{G}_t $
  such that
   $$(\phi - \mathrm{ad}D)(\Gamma_k-\Gamma_{1'})=0
   \; \mbox{for all}\  k\in Y_1\setminus 1'.$$
  \end{lemma}
\begin{proof}  By Lemma  \ref{gt3.2.1} (i) and (iii) it suffices to consider
the setting $t<n-1.$
      By Lemma \ref{lemma2.1}, one may assume that
      \begin{equation}\phi(\Gamma_k-\Gamma_{1'})=\sum_{r\in Y_0}f_{rk} D_r
       \;\mbox{where}\ k\in Y_1\setminus 1';\ f_{rk}\in \Lambda(n).
       \label{ge3.2.13}
       \end{equation}
Write
      \begin{equation}f_{rk}=\sum_{|u|=t+1}c_{u,r,k}x^u\quad \mbox{where}\
       c_{u,r,k}\in \mathbb{F}.\label{ge3.2.15}
       \end{equation}
     Discussing just as in the proof of \cite[Lemma 4.2.4]{g6}, we may obtain that
      $$\sum_{|u|=t+1}(\delta_{k\in u}-\delta_{1'\in u})c_{u,r,l}x^u=
      \sum_{|u|=t+1}(\delta_{l\in u}-\delta_{1'\in u})c_{u,r,k}x^u.$$
      Since $\{x^u\mid u\in \mathbb{B}\}$ is an $\mathbb{F}$-basis
      of $\Lambda(n),$ it follows
       that
      \begin{equation}(\delta_{k\in u}-\delta_{1'\in u})c_{u,r,l}=
      (\delta_{l\in u}-\delta_{1'\in u})c_{u,r,k}
      \quad \mbox{for}\ r\in Y_0,\ k,l\in Y_1\setminus
      1'.\label{ge3.2.16}
      \end{equation}
      Suppose that $ c_{u,r,k} $ is any a nonzero coefficient in (\ref{ge3.2.15}),
      where $ |u|=t+1< n,\, r\in Y_0 $ and $ k\in Y_1\setminus 1^{'}.$
      Note that $ |u|\geq 1.$
We proceed  in two steps to show that $ \delta_{k\in
u}+\delta_{1^{'}\in u}=1.$\\

 \noindent\textit{Case }(i): $ |u|\geq 2.$  If $ 1^{'}\not\in u,$ one may find
  $ l\in u\setminus k.$  Then (\ref{ge3.2.16}) shows that $ \delta_{k\in u}=1;$ that
 is,  $ k\in u.$ If $ 1^{'}\in u,$ noting that $ |u|\leq n-1,$ one may find
$  l\in Y_1 \setminus u.$  Then (\ref{ge3.2.16}) shows that $
\delta_{k\in u}=0; $ that is, $ k\not\in u.$  Summarizing, for any
nonzero coefficient $ c_{u,r,k} $ in (\ref{ge3.2.15}), we have $
\delta_{k\in u}+\delta_{1^{'}\in u}=1.$\\

 \noindent\textit{Case }(ii): $ |u|=1.$ Since $ |u|=1,$ the case of $ k\in u$ and $1^{'}\in
 u$ does not occur. If $ k\not\in u $ and
 $ 1^{'}\not\in u,$ then there is $ l\in u,$ since
  $ |u|=1.$ Then by (\ref{ge3.2.16}), we get $ c_{u,r,k}=0,$ this is a
 contradiction. Hence, we have $ \delta_{k\in u}+\delta_{1^{'}\in
 u}=1.$

Then, just like in the proof of \cite[Lemma 4.2.4]{g6}, we can rewrite
(\ref{ge3.2.15}) as follows
    $$f_{rk}=\sum_{ 1'\in u,k\not\in u} c_{u,r,k}x^u+
     \sum_{ 1'\not\in u, k\in u} c_{u,r,k}x^u.$$
   Now, following the corresponding part of the proof for \cite[Lemma 4.2.4]{g6},
   one may find $D\in \mathcal{G}_{t}$ such that
   $(\phi-\mathrm{ad}D)(\Gamma_k-\Gamma_{1'})=0
   \; \mbox{for all}\ k \in Y_1\setminus 1'.$ The proof is complete.
   \end{proof}

    Let us consider the case of odd $ \mathbb{Z}$-degree.

\begin{lemma}\label{gt3.2.5}

Let $\phi \in \mathrm{Der}_{t}(\mathcal{S},W_{\bar{1}})
   $ where $t> 0$ is odd.  If $\phi(\mathcal{S}_{-1})=0,$
  then there is $D\in \mathcal{G}_t $
   such that
   \begin{equation}(\phi - \mathrm{ad}D)(\Gamma_1+\Gamma_k)=0
   \; \mbox{for all}\  k\in Y_1.\label{ge3.2.21}
   \end{equation}
\end{lemma}
\begin{proof} Deleting the part (ii) in the proof of \cite[Lemma 4.2.5]{g6}, we
obtain our proof.
\end{proof}

 For our purpose, we need still the following three reduction
 lemmas.

\begin{lemma}\label{gt3.2.6} 

 Suppose that $\phi \in
\mathrm{Der}_{t}(\mathcal{S},W_{\bar{1}})
   $ and $\phi(\mathcal{S}_{-1})=0, $ where $t> 0$ is odd.
   If $\phi(\Gamma_1+\Gamma_k)=0 $
   for all $k\in Y_1,$
   then $\phi(\mathcal{S}_{0})=0.$
\end{lemma}
\begin{proof} Following parts (i) and (ii) in the proof
of \cite[lemma 4.2.6]{g6},
one may show that
   $\phi(\Gamma_1-\Gamma_i)=0$
and
    $\phi(x_iD_j)=0 $
   for all $i,j\in Y_0$ with $i\neq j.$

   To show that $\phi(x_kD_l)=0$ for $k,l\in Y_1$ with $k\neq l,$ just as in the part (iii) of
   the proof of \cite[Lemma 4.2.6]{g6},
   it suffices to consider separately two cases $\mathrm{zd}(\phi)=1 $
   and $\mathrm{zd}(\phi)\geq 3. $
Now Lemma \ref{gt2.1.7} ensures that $\phi(\mathcal{S}_{0})=0.$
\end{proof}

Analogous to \cite[Lemma 4.2.7]{g6}, one may prove the following

 \begin{lemma}\label{gt3.2.7}

 Suppose that $\phi \in
   \mathrm{Der} (\mathcal{S},W_{\bar{1}})
   $ is a homogeneous derivation of nonnegative even $\mathbb{Z}$-degree
  and $\phi(\mathcal{S}_{-1}+\mathcal{T_{\mathcal{S}}})=0.$
   Then $\phi(\mathcal{S}_{0})=0.$
 \end{lemma}

Now we are able to characterize the homogeneous derivation
 space of nonnegative $\mathbb{Z}$-degree. Using
 Lemmas \ref{gt3.2.4}--\ref{gt3.2.7}, Corollaries \ref{gt3.1.4} and  \ref{gt3.2.3}, and
  Proposition \ref{lemma2.3}, one may prove the following result
  (cf. \cite[Proposition 4.2.9]{g6}).

\begin{proposition}\label{gt3.2.8}

 $ \mathrm{Der}_t(\mathcal{S},W_{\bar{1}})= \mathrm{ad}(W_{\bar{1}})_t \; for \; t\geq
0.$
\end{proposition}

 As an application of Proposition \ref{gt3.2.8}, we have:

 \begin{proposition}\label{gt3.2.9}

$ \mathrm{Der}_t(\mathcal{S},S_{\bar{1}})=
\mathrm{ad}(\overline{S}_{\bar{1}})_t \; for \; t\geq 0. $
\end{proposition}

 \begin{proof} Since
 $ S$ is an ideal of $\overline{S}$ (see \cite[p. 139]{ZZ1}),
 $\mathrm{ad}(\overline{S}_{\bar{1}})_t
\subset\mathrm{Der} _{t}(\mathcal{S},S_{\bar{1}}).$
 Let $\phi\in\mathrm{Der} _{t}(\mathcal{S}, S_{\bar{1}}).$ View $\phi$ as
 a derivation of $\mathrm{Der} _{t}(\mathcal{S}, W_{\bar{1}}).$
 Then by Proposition \ref{gt3.2.8}, there is $D\in (W_{\bar{1}})_t$ such that
 $\phi=\mathrm{ad}D\in \mathrm{Der} _{t}(\mathcal{S}, W_{\bar{1}}),$
 and therefore,
 $\phi=\mathrm{ad}D\in \mathrm{Der} _{t}(\mathcal{S}, S_{\bar{1}}).$ Let
 $\mathrm{Nor}_{W_{\bar{1}}}(\mathcal{S}, S_{\bar{1}}):=\{x\in W_{\bar{1}}\; |\;
 [x, \mathcal{S}]\subset S_{\bar{1}} \}.$
 Clearly,
 $D\in \mathrm{Nor}_{W_{\bar{1}}}(\mathcal{S}, S_{\bar{1}})_{t}.$ Therefore,
 it suffices to show that
 $\mathrm{Nor}_{W_{\bar{1}}}(\mathcal{S}, S_{\bar{1}})_{t}
 \subset (\overline{S}_{\bar{1}})_{t}.$
  Let $E$ be an arbitrary element of
  $\mathrm{Nor}_{W_{\bar{1}}}(\mathcal{S}, S_{\bar{1}})_{t}.$
  If $t=0$ then $E\in W_{\overline{1}}\cap W_{0},$ which implies
  that $\mathrm{div}(E)=0 $
  and therefore, $E\in \overline{S}_{\overline{1}}.$ Now suppose that $t>0.$
  Note that
  $\mathrm{div} ([E, \mathcal{S}_{-1}] )=0. $
  It follows
    that $D_i(\mathrm{div}(E))=0 $ for all $i\in Y_0.$ This implies that
  $\mathrm{div}(E)\in \Lambda(n)_{t}.$ Similarly, $[E,\mathcal{S}]\subset
  S_{\overline{1}}$ implies that $[\mathrm{div}(E),
  \mathcal{S}]=0.$ In particular, $[\mathrm{div}(E),
  \mathcal{T}_{\mathcal{S}}]=0.$ Since $\mathrm{div}(E)\in
  \Lambda(n)_{t} ,$  one gets $[\mathrm{div}(E),
  \mathcal{T}]=0. $   Keeping in mind  $\mathrm{div}(E)\in \Lambda(n)_{\overline{1}},$  one
  may easily deduce that $\mathrm{div}(E)=0 $ (c.f \cite[Proposition  4.2.10]{g6}).
\end{proof}

  In the following we first determine the homogeneous
derivations of negative $\mathbb{Z}$-degree from $\mathcal{S}$
into $ W_{\bar{1}}.$ This combining with Proposition \ref{gt3.2.8} will give the structure of the derivation space
$\mathrm{Der}(\mathcal{S} ,W_{\bar{1}}).$
 The following lemma tells us that a $\mathbb{Z}$-degree $-1$ derivation
  from $\mathcal{S} $ into $ W_{\bar{1}} $ is completely
 determined by its action on $ \mathcal{S}_0.$

 \begin{lemma}\label{gt3.3.1}

Suppose that $ \phi\in \mathrm{Der}_{-1}(\mathcal{S},
W_{\bar{1}})$ and that $\phi(\mathcal{S}_0)=0.$ Then $ \phi=0.$
   \end{lemma}
\begin{proof}
We first show that $ \phi(\mathcal{R})=0.$ By the definition of  $
D_{il} ,$
\begin{equation}
D_{il}(x^{(2\varepsilon_i)}x_k)=x_ix_kD_l+\delta_{kl}x^{(2\varepsilon_i)}D_i
\quad \mbox {for all} \; i\in Y_0,\,k,l\in Y_1. \label{ge3.3.1}
\end{equation}
 We shall use the
 following simple fact (by Lemma \ref{lemma2.1}):
 $$ \phi(\mathcal{S}_1)\subseteq \mathcal{G}_0.$$
 We may assume that
 \begin{equation} \phi(D_{il}(x^{(2\varepsilon_i)}x_k))=\sum_{k\in Y_1, r\in
 Y_0}c_{k,r}x_kD_r. \label{ge3.3.2}
 \end{equation}
Given  $ i\in Y_0,\, k,l\in Y_1, $ take $ j\in Y_0\setminus i.$
Then
\begin{equation*} [\Gamma_i-\Gamma_j,
D_{il}(x^{(2\varepsilon_i)}x_k)]=D_{il}(x^{(2\varepsilon_i)}x_k).
\end{equation*}
Applying $ \phi $ to the equation above and then combining that
with (\ref{ge3.3.2}), one may obtain by a comparison of
coefficients that
$$ c_{k,r}=0 \quad \mbox{for} \; k\in Y_1,\; r\in Y_0\setminus j. $$
Hence, by (\ref{ge3.3.2}), we may obtain
\begin{equation}\phi(D_{il}(x^{(2\varepsilon_i)}x_k))=\sum_{k\in Y_1}c_{k,j}x_kD_j.
\label{ge3.3.4}
\end{equation}
\noindent  \textit{Case} (i): $ k\neq l.$ Then by (\ref{ge3.3.1}),
we have $D_{il}(x^{(2\varepsilon_i)}x_k)=x_ix_kD_l $ and
\begin{equation} [\Gamma_i+\Gamma_k,
D_{il}(x^{(2\varepsilon_i)}x_k)]=2D_{il}(x^{(2\varepsilon_i)}x_k).
\label{ge3.3.5}
\end{equation}
 Applying $ \phi $ to (\ref{ge3.3.5}) and using (\ref{ge3.3.4}), one may obtain by comparing coefficients that
 $$ c_{s,j}=0 \quad \mbox{for} \; s\in Y_1.$$
Consequently, $ \phi(D_{il}(x^{(2\varepsilon_i)}x_k))=0. $\\

\noindent  \textit{Case} (ii): $ k=l.$ Then by (\ref{ge3.3.1}), we
have
$$ D_{ik}(x^{(2\varepsilon_i)}x_k)=x^{(2\varepsilon_i)}D_i+x_ix_kD_k. $$
and $$ [\Gamma_i+\Gamma_k
,D_{ik}(x^{(2\varepsilon_i)}x_k)]=D_{ik}(x^{(2\varepsilon_i)}x_k).$$
Applying $ \phi $ to the equation above and using (\ref{ge3.3.4}),
one may obtain by comparing coefficients that
$$ c_{s,j}=0 \quad \mbox{for} \; s\in Y_1\setminus k. $$
Hence, By (\ref{ge3.3.4}), we may obtain
\begin{equation}\phi(D_{ik}(x^{(2\varepsilon_i)}x_k))=c_{k,j}x_kD_j. \label{ge3.3.6}
\end{equation} For $ i\in Y_0,\, k,l\in Y_1,$ choose $ q\in
Y_1\setminus k.$ Then
 $ m\Gamma_q+\Gamma^{''}\in \mathcal{S}_0 $ and
$$ [m\Gamma_q+\Gamma^{''}, D_{ik}(x^{(2\varepsilon_i)}x_k)]=D_{ik}(x^{(2\varepsilon_i)}x_k).$$
Applying $ \phi $ and using (\ref{ge3.3.6}), one gets
$$\phi(D_{ik}(x^{(2\varepsilon_i)}x_k))= [m\Gamma_q+\Gamma^{''},\phi(D_{ik}(x^{(2\varepsilon_i)}x_k))]=
-\phi(D_{ik}(x^{(2\varepsilon_i)}x_k)).$$ Consequently, $
\phi(D_{ik}(x^{(2\varepsilon_i)}x_k))=0 \;\mbox {for all}\; i\in
Y_0, k\in Y_1.$

 We want to show that $ \phi(D_{ij}(x_ix_kx_l))=0$  for $i,j\in Y_0,
 k,l\in Y_1.$
 By a same argument, we can also obtain that
 \begin{equation} \phi(D_{ij}(x_ix_kx_l))=c_{kj}x_kD_j \quad \mbox{for}\; k\in Y_1, j\in Y_0\setminus i.
 \label{ge3.3.7}
 \end{equation}
 Note that $ n\Gamma_q+\Gamma^{'}\in \mathcal{S}_0 $ for $ q\in Y_0$ and
 that
\begin{equation} [n\Gamma_q+\Gamma^{'}
,D_{ij}(x_ix_kx_l)]=(2-n\delta_{q,j})D_{ij}(x_ix_kx_l)\quad
\mbox{for} \; i,j\in Y_0,k,l\in Y_1. \label{ge3.3.8}
\end{equation}
Applying $ \phi $ to (\ref{ge3.3.8}) and using (\ref{ge3.3.7}),
one gets
$$ c_{kj}=0 \quad \mbox{for}\; k\in Y_1,j\in Y_0\setminus i. $$
Consequently, $ \phi(D_{ij}(x_ix_kx_l))=0.$

It remains to show that $\phi(\mathcal{Q})=0.$ But this can be verified completely
analogous to the proof of \cite[Lemma 4.3.1]{g6}.
   \end{proof}

 Using Lemma \ref{gt3.3.1} we can determine the derivations of
 $\mathbb{Z}$-degree $-$1.

 \begin{proposition}\label{gt3.3.2}

$
\mathrm{Der}_{-1}(\mathcal{S},W_{\bar{1}})=\mathrm{ad}(W_{\bar{1}})_{-1}.$
In particular, $
\mathrm{Der}_{-1}(\mathcal{S},S_{\bar{1}})=\mathrm{ad}(\overline{S}_{\bar{1}})_{-1}.$
 \end{proposition}
 \begin{proof}
Let $\phi\in \mathrm{Der}_{-1}(\mathcal{S},W_{\bar{1}}).$  Let $
k\in Y_1.$
 Assume that
 \begin{equation} \phi(\Gamma_k+\Gamma_1)=\sum_{r\in Y_1}c_{kr}D_r\quad
 \mbox{where}\ c_{kr}\in
 \mathbb{F}. \label{ge3.3.15}
 \end{equation}
 Let $l\in Y_1\setminus k.$ Then
 $[\Gamma_k+\Gamma_1
,\phi(\Gamma_l+\Gamma_1)]=[\Gamma_l+\Gamma_1
,\phi(\Gamma_k+\Gamma_1)].$
 By (\ref{ge3.3.15}), $ c_{kl}=0 $ whenever $k,l\in Y_1\;\mbox{with} \; k\neq l. $ It follows
 that
$ \phi(\Gamma_k+\Gamma_1)=c_{kk}D_k \quad \mbox{where} \; c_{kk}\in
 \mathbb{F}. $
 Obviously, $ [\Gamma_k+\Gamma_1, x_kD_l]=x_kD_l $  for
 $ k,l\in Y_1 $ with $ k\neq l. $
 Then
 $$c_{kk}D_l+[\Gamma_k+\Gamma_1, \phi(x_kD_l)]=\phi(x_kD_l).$$
 Since $ \phi(x_kD_l)\in (W_{\bar{1}})_{-1},$
it follows that $\phi(x_kD_l)=c_{kk}D_l.$ Put
$ \psi:=\phi-\sum_{r\in Y_1}c_{rr}\mathrm{ad}D_r.$ Then
 \begin{equation} \psi(\Gamma_k+\Gamma_1)
 =\psi(x_kD_l)=0 \quad \mbox{for} \; k,l\in Y_1, \,
  k\neq l. \label{ge3.3.17}
  \end{equation}
 We next show that
\begin{equation} \psi(x_iD_j)=0 \quad \mbox{for}
 \, i,j\in Y_0, \; i\neq j. \label{ge3.3.18}
\end{equation}
 Take $ r\in Y_0\setminus {\{i,j\}}.$ Then
$[\psi(\Gamma_r+\Gamma_q),x_iD_j]+[\Gamma_r+\Gamma_q,\psi(x_iD_j)]=0. $
 Since $ \psi(\Gamma_r+\Gamma_q)\in (W_{\bar{1}})_{-1},$ we have $[\psi(\Gamma_r+\Gamma_q)
,x_iD_j]=0.$ Consequently, $ [\Gamma_q
,\psi(x_iD_j)]=0 $
 for all $ q\in Y_1.$  Hence $ \psi(x_iD_j)=0,$ since  $ \psi(x_iD_j)\in
 (W_{\bar{1}})_{-1}.$

 In the same way we can verify that
\begin{equation} \psi(\Gamma_1-\Gamma_j)=0 \quad \mbox {for all}\; j\in
Y_0\setminus 1. \label{ge3.3.20}
\end{equation}
By (\ref{ge3.3.17})--(\ref{ge3.3.20}), we have
$\psi(\mathcal{S}_0)=0.$ It follows from Lemma \ref{gt3.3.1}  that $
\psi=0 $ and hence $ \phi\in \mathrm{ad}(W_{\bar{1}})_{-1}.$  This completes the proof.
\end{proof}

To compute the derivation  of $\mathbb{Z}$-degree less than $-1$
from $\mathcal{S}$ into $W_{\bar{1}},$ we establish the following
lemma.
 \begin{lemma}\label{gt3.3.3}

Suppose that $ \phi\in \mathrm{Der}_{-t}(\mathcal{S},W_{\bar{1}})$
with $ t>1$ and that
 $ \phi(D_{ij}(x^{((t+1)\varepsilon_i)}))=0 $ for all $ i,j\in Y_0.$
Then $ \phi=0.$
 \end{lemma}
  \begin{proof}
First claim that $ \phi(\mathcal{Q})=0.$ To that aim, we proceed
by induction on $ q $ to show that
\begin{equation}
\phi(D_{ij}(x^{(q\varepsilon_i)}))=0\quad \mbox {for all}\;  i,j\in
Y_0 \; \mbox{with}\ i\neq j. \label{ge3.3.21}
\end{equation}
Without loss of generality  suppose that $ q>t+1
$ in the following. By inductive hypothesis and Lemma
\ref{lemma2.1},
 $\phi(D_{ij}(x^{(q\varepsilon_i)}))\in \mathcal{G}_{q-t-2}.$ Thus one
may write
 \begin{equation}
\phi(D_{ij}(x^{(q\varepsilon_i)}))=\sum_{r\in
Y,|u|=q-t-1}c_{u,r}x^{u}D_r \quad \mbox{where} \; c_{u,r}\in
\mathbb{F}.\label{ge3.3.22}
\end{equation}
If $ q-t\geq 3,$ proceeding just as Case (i) in the proof of \cite[Lemma 4.3.3]{g6},
one may show that
$ \phi(D_{ij}(x^{(q\varepsilon_i)}))=0.$
Suppose that $ q-t<3.$ Note that $ q>t+1.$ Then
rewrite (\ref{ge3.3.22}) as
\begin{equation}
\phi(D_{ij}(x^{(q\varepsilon_i)}))=\sum_{l\in Y_1, r\in
Y_0}c_{l,r}x_lD_r \quad \mbox{where} \; c_{l,r}\in \mathbb{F}.
\label{ge3.3.23}
\end{equation}
For any fixed coefficient $ c_{l_0,r_0}$ in (\ref{ge3.3.23}).
Choose $ k\in Y_1\setminus l_0 ,\ s\in Y_0\setminus {\{r_0,i\}},$
since $n,m\geq 3.$

 If $ s=j,$ then $
[\Gamma_s+\Gamma_k ,
D_{ij}(x^{(q\varepsilon_i)})]=-D_{ij}(x^{(q\varepsilon_i)}).$
Applying $ \phi $ to the equation above and then combining that
with (\ref{ge3.3.23}), one may obtain by a comparison of
coefficients of $ x_{l_0}D_{r_0} $ that
$$ c_{l_0,r_0}=0 \quad \mbox{for}\; l_0\in Y_1, r_0\in Y_0. $$
Consequently, $ \phi(D_{ij}(x^{(q\varepsilon_i)}))=0. $

 If $ s\neq j
,$ then $ [n\Gamma_s+\Gamma^{'} , D_{ij}(x^{(q\varepsilon_i)})]=0.
$ Applying $ \phi $ and then combining that with (\ref{ge3.3.23}),
one may obtain by a comparison of coefficients of $
x_{l_0}D_{r_0}$ that
$$ c_{l_0r_0}=0 \quad \mbox{for}\; l_0\in Y_1, r_0\in Y_0. $$
Consequently, $ \phi(D_{ij}(x^{(q\varepsilon_i)}))=0. $ Thus
(\ref{ge3.3.21}) holds for all q and therefore, $
\phi(\mathcal{Q})=0.$

 We next prove that $ \phi(\mathcal{R})=0.$ Since
$ \mathcal{R}\subseteq \mathcal{S}_1,$ $ \mathrm{zd}(\phi)\leq
-2,$ it
 suffices to consider the case that $ \mathrm{zd}(\phi)=-2.$ Note
 that $ \phi(\mathcal{S}_1)\subset \mathcal{S}_{-1}.$
 For $
k,l\in Y_1, i\in Y_0,$ take $ q\in Y_0 \setminus {\{i,j\}}.$
Then $ n\Gamma_q+\Gamma^{'}\in \mathcal{S}_0 $ and
$ [n\Gamma_q+\Gamma^{'},D_{il}(x^{(2\varepsilon_i)}x_k)]=0.$
Since $ \phi(D_{il}(x^{(2\varepsilon_i)}x_k))\in
(W_{\bar{1}})_{-1},$ it follows that
$$ \phi(D_{il}(x^{(2\varepsilon_i)}x_k))=-[n\Gamma_q+\Gamma^{'}
,\phi(D_{il}(x^{(2\varepsilon_i)}x_k))]=0\quad \mbox {for all} \;
i\in Y_0,\; k,l\in Y_1.$$
  Obviously, $$[n\Gamma_q+\Gamma^{'}
,D_{ij}(x_ix_kx_l)]=2D_{ij}(x_ix_kx_l).$$ Applying $ \phi ,$ one
gets,
$$ 2\phi(D_{ij}(x_ix_kx_l))=[n\Gamma_q+\Gamma^{'}
,\phi(D_{ij}(x_ix_kx_l))]=-\phi(D_{ij}(x_ix_kx_l)),$$ since $
\phi(D_{ij}(x_ix_kx_l))\in (W_{\bar{1}})_{-1}. $ The assumption $p\not=3$ ensures that $
\phi(D_{ij}(x_ix_kx_l))=0. $ By Lemma \ref{gt2.1.8}, $ \phi=0.$
\end{proof}

Finally, we are to determine the homogeneous derivations of
$\mathbb{Z}$-degree less than $-1$ from $\mathcal{S}$ into
$W_{\bar{1}}.$
 \begin{proposition}\label{gt3.3.4}

  $\mathrm{Der}_{-t}(\mathcal{S},W_{\bar{1}})=0\; \mbox{for} \;t>1.$  In
particular, $ \mathrm{Der}_{-t}(\mathcal{S},S_{\bar{1}})=0 \;
\mbox{for} \; t>1.$
\end{proposition}
\begin{proof} Let $ \phi\in
\mathrm{Der}_{-t}(\mathcal{S},W_{\bar{1}}).$ Assert  that
 $$ \phi(D_{ij}(x^{((t+1)\varepsilon_i)}))=0 \quad \mbox {for
 all}\; i,j\in Y_0.  $$
 Recall
$ \Gamma^{'}=\sum_{r\in Y_1}\Gamma_r.$ Choose $ q\in Y_0\setminus
{\{i,j\}},$ since $ m\geq 3.$ Clearly, $ n\Gamma_q+\Gamma^{'}\in
\mathcal{S}_0.$
 Then $$[ n\Gamma_q+\Gamma^{'},D_{ij}(x^{((t+1)\varepsilon_i)})]=0.$$
Applying $ \phi,$ one gets
 $$ 0=[ n\Gamma_q+\Gamma^{'}, \phi(D_{ij}(x^{((t+1)\varepsilon_i)}))]=-\phi(D_{ij}(x^{((t+1)\varepsilon_i)})),$$
since $\phi(D_{ij}(x^{((t+1)\varepsilon_i)}))\in
(W_{\bar{1}})_{-1}.$ Consequently, $
\phi(D_{ij}(x^{((t+1)\varepsilon_i)}))=0.$
By Lemma \ref{gt3.3.3}, $ \phi=0.$ The proof is complete.\end{proof}

Now we can describe the derivation spaces
$\mathrm{Der}(\mathcal{S},W_{\bar{1}})$ and
$\mathrm{Der}(\mathcal{S},S_{\bar{1 }}).$
\begin{theorem}\label{gt3.3.5}

 $\mathrm{Der}(\mathcal{S},W_{\bar{1}} )
 =\mathrm{ad}W_{\bar{1}}.$
\end{theorem}

 \begin{proof} This is  a direct consequence of
 Propositions \ref{gt3.2.8}, \ref{gt3.3.2}, \ref{gt3.3.4}. \end{proof}

\begin{theorem}\label{gt3.3.6}

$\mathrm{Der}(\mathcal{S},S_{\bar{1}})
 =\mathrm{ad}\overline{S}_{\bar{1}}.$
\end{theorem}

 \begin{proof} This is  a direct consequence of
 Propositions \ref{gt3.2.9}, \ref{gt3.3.2}, \ref{gt3.3.4}. \end{proof}

\end{document}